\documentclass[review,11.5pt]{article}

\usepackage{lineno,hyperref}
\modulolinenumbers[5]

\usepackage{amssymb}
\usepackage[fleqn]{amsmath}
\usepackage[margin=1in]{geometry}
\usepackage{multirow}
\usepackage[table,xcdraw]{xcolor}
\usepackage{wrapfig}
\usepackage{float}
\usepackage{caption}
\usepackage{subcaption}
\usepackage{authblk}
\usepackage{booktabs}
\usepackage{tablefootnote}
\usepackage{threeparttable}
\usepackage{lineno}
\usepackage{natbib}
\usepackage{graphicx}
\usepackage{placeins}

\title{A Multi-Objective Capacity-Constrained Optimization of Corn Planting Scheduling}

\author[1]{Mingshi Cui}

\author[2]{Kunting Qi}

\author[3]{Byran Smucker\thanks{corresponding author, current email and affiliations: bsmucke1@hfhs.org; Henry Ford Health + Michigan State University Health Sciences, Detroit, MI, USA; Dept of Public Health Sciences, Henry Ford Health, Detroit, MI, USA; Dept of Epidemiology \& Biostatistics, College of Human Medicine, Michigan State University, East Lansing, MI, USA}}

\author[4]{Durai Sundarmoorthi}

\affil[1]{Rutgers University, Department of Statistics, New Brunswick, NJ USA}

\affil[2]{University of Illinois Chicago, Department of Computer Science, Chicago, IL USA}

\affil[3]{Miami University, Department of Statistics, Oxford, OH USA}

\affil[4]{Washington University in St. Louis, Olin Business School, St. Louis, MO USA}

\begin{document}

\maketitle

\begin{abstract}
This article describes an improved set of solutions to the problems presented in the 2021 Syngenta Crop Challenge in Analytics \citep{Syngenta2021}. In particular, we use multiobjective optimization and predictive modeling methods to determine a corn planting schedule. The problem involves the following objectives: i. minimize the median and maximum absolute difference between weekly harvest quantity and the storage capacity, the number of nonzero harvest weeks, and the total amount of corn wasted. This is accomplished while respecting planting windows, expected harvest amounts, the growing degree units  required to bring seeds to harvest, and historical weather data. We used a Long Short-Term Memory model to predict growing degree units for 2020 and 2021, based on historical data. Then, we used a genetic algorithm, and an extensive search of the tuning parameter space, to produce a Pareto front of solutions for three distinct optimization models related to the Challenge. We evaluate the quality of the Pareto fronts for each model, and use the results to choose a preferred model and final solution. We also provide comparisons between our final solutions, previous solutions submitted to the Challenge, and solutions from other groups.
\end{abstract}

\noindent%
{\it Scheduling; Agricultural Optimization; Hypervolume; Genetic Algorithms}


\section{Introduction and Motivation} \label{sec:intro}

Corn is a key food commodity that provides benefits to human beings both directly through its consumption and indirectly through animal consumption. As commercial corn-growing technology has improved, the obvious benefits of food and energy availability are limited by challenges regarding storage capacity. To tackle this problem, Syngenta - a leading agriculture company, posed an analytics challenge to the operations research community in their 2021 Syngenta Crop Challenge \citep{Syngenta2021}. Syngenta and other similar agriculture companies grow and harvest corn in South America multiple times within a year by taking advantage of the longer growing season there. Such an advantage requires careful planning of growing, harvesting, and storage operations. The modeling framework in this article reflects growing, harvesting, and storage practices of companies like Syngenta. Given a specified storage capacity, along with historical data about weather, potential planting intervals for different varieties of corn, the Growing Degree Units (GDUs) needed to bring the planted seeds to harvest, and the amount of harvest expected, the main question of interest is: when should a set of seed populations be planted in order to minimize waste and provide as consistent a weekly harvest quantity as possible? More specifically, we wish to minimize (a) the median difference between weekly harvest quantity and the capacity across all weeks; (b) the maximum difference between the weekly harvest quantity and the capacity across all weeks; and (c) the number of weeks that a non-zero harvest is realized. Though it was not a part of the original challenge, we have added a fourth criterion to minimize: (d) the total amount of corn wasted because the weekly harvest exceeds capacity. These criteria should be simultaneously optimized while respecting the assigned planting window for each seed population, expected harvest, GDUs required, and historical weather data. This is the first Scenario presented in the Challenge. The second Scenario generalizes the problem to require a reasonable capacity specification instead of taking it as given.

\begin{figure}[hbt!]
    \centering
    \includegraphics[width=0.95\textwidth]{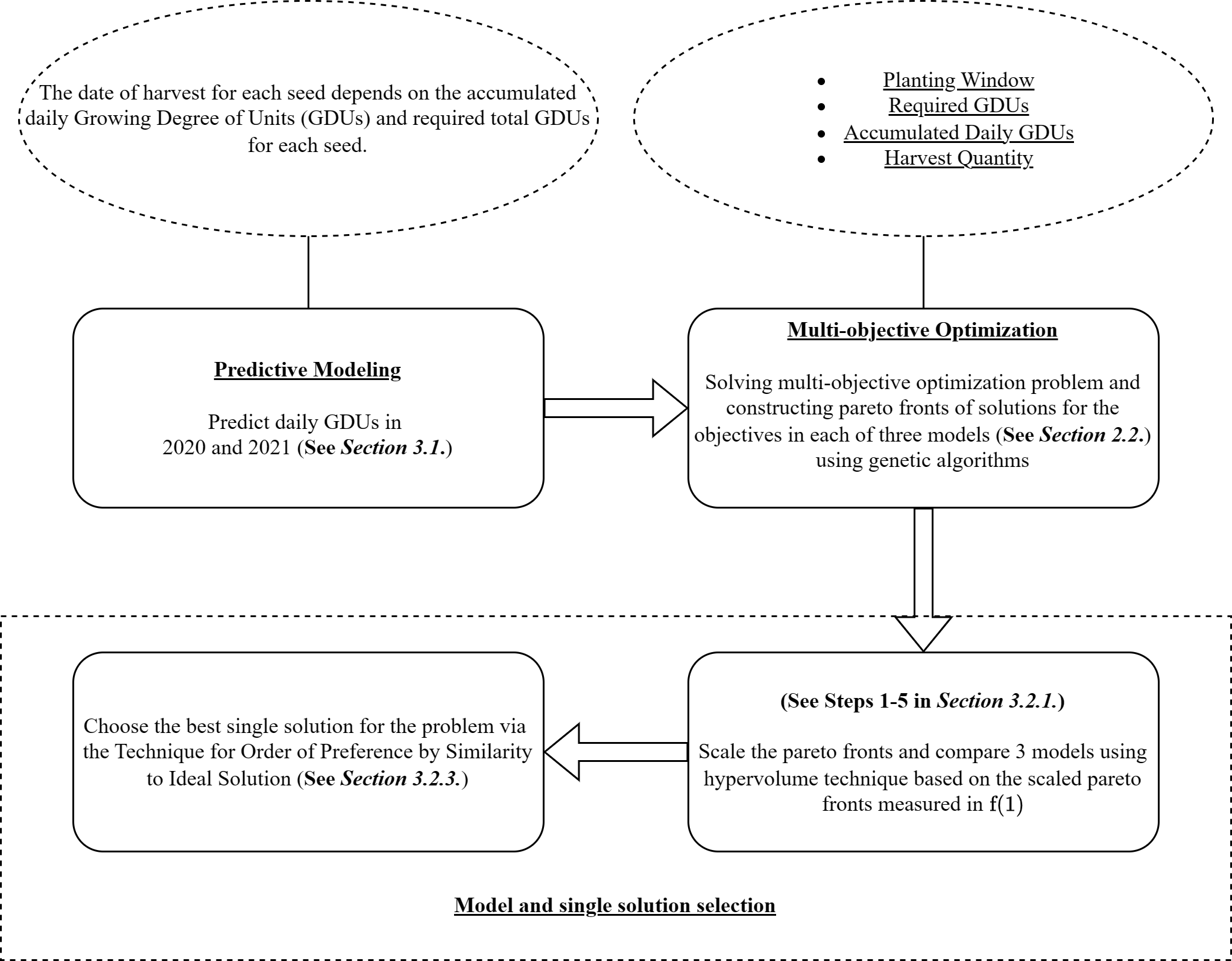}
    \caption{Overview of our solution methodology.}
    \label{fig:workflow}
\end{figure}

In this article, we formulate these problems as novel applications of multiobjective optimization, which we solve using a genetic algorithm. The problems have several challenging and unique aspects. First, in order to obtain a solution, we must use historical data to predict the Growing Degree Units (GDUs) over the planting year. This adds a predictive modeling aspect to the optimization problem, and we consider both a simple averaging model as well as a neural network-based approach. Second, we sharpen the basic optimization model by considering two related models that reflect the asymmetry inherent in the median and maximum difference objectives. We then construct Pareto fronts based on all three models but assess these Pareto fronts based on the original model, since it most explicitly reflects the problem outlined in the Challenge. Interestingly, we find that we always obtain better results using one of the related models, instead of using the original model. Third, on the multiobjective algorithm side, we perform an extensive round of parameter tuning to improve the optimization results, and use the hypervolume indicator to discriminate between competing Pareto fronts, before using another method from the literature to select a final solution. See Figure \ref{fig:workflow} for an overview of our solution strategy. 
We compare our solutions with those that we originally submitted as part of the second-place entry to the Challenge, as well as with results from two other teams. 

Our efforts are situated within the optimization literature as applied to agriculture; in particular, work which uses several elements of analytics, along with traditional optimization, to guide decision-making \citep[see][for a review]{LoweTimothyPreckelPaul}. We highlight a few particularly relevant works. 
\citet{JonesPhilipLoweTimothy} used a two-period sequential production scheduling problem for an agribusiness which has operations in North America and South America. 
In the first period, they have production operations in North America between April and September. In the second period, they produce in South America from October to March. This is reminiscent of the present problem which has two separate grow sites. 
\citet{bansal2017} proposed an optimization-based approach to estimate parameters for yield distribution based on expert judgement about the potential yield of different seed corn varieties. 
\citet{bansal2017} also indicate that a lack of data is a limiting factor in this area of research. This situation has been improved after Syngenta – a leading agribusiness - won the 2015 Edelman Award given by the Institute for Operations Research and the Management Sciences (INFORMS) for \citet{ByrumJosephDavis}, which used simulation and optimization in seed development and the production process. Since then, Syngenta organized a series of analytics challenges to solve relevant business problems, providing data and context that has stimulated substantial work. 
Such work includes \citet{Ansarifar2020} and \citet{sundaramoorthi2022machine}, which are characterized by the need to combine analytics strategies---including machine learning and optimization---to solve complex problems. Our work carries on that tradition. In particular, we elaborate upon the second-place entry of the 2021 Syngenta Crop Challenge, and presents a novel solution to a problem which to our knowledge has only recently been considered in the literature by two other  participating teams \citep{sajid2022optimizing, lizhiwang}. 

The rest of the article continues as follows. Section \ref{sec:problem_data} provides a detailed description of the problem as well as the mathematical formulation of the three optimization models that we consider. Section \ref{sec:methods} describes the methodological approach we took to solve the problem, including a description of the predictive model we used, the optimization procedure, and how we tuned the algorithm. In Section \ref{sec:results}, we demonstrate our methodology and provide particular solutions. Finally, we close with a short discussion in Section \ref{sec:conclusion}, including comparison with solutions from other published work.

\section{Problem Description and Formulation} \label{sec:problem_data}

In this section we provide a general description of the problem and data that was provided, and then a mathematical formulation of the optimization problems we are trying to solve.

\subsection{Description}

The original motivation for this work was the 2021 Syngenta Crop Challenge. The problem was to provide a corn planting schedule that results in a weekly harvest quantity consistently close to capacity. Two datasets were provided. First, information about a set of $2{,}569$ seed populations, including the site at which each population was to be planted (site 0 or 1); a window of time within which each population could be planted (early planting dates ranged from January 1, 2020 to December 29, 2020; late planting dates ranged from January 22, 2020 to February 16, 2021); the number of required GDUs needed to bring the seeds to harvest (range from 649 to 1{,}414); and the number of ears that the seeds would produce once they received the necessary GDUs, in Scenario 1 or Scenario 2. 
Second, Syngenta provided the historical GDUs for each day from 2009 to 2019, at sites 0 and 1. The data provided is summarized in Table \ref{tab:data_description}. Thus, each seed population must be planted at its specified site, but the seed population can be planted under either Scenario 1 or Scenario 2, which have differing harvest quantities. Furthermore, the date when a seed population is ready to harvest is determined by the site-specific daily GDUs and its own required GDUs. Note that the raw data provided to Challenge participants is proprietary and 
cannot be shared publicly.


\begin{table}[h]
\centering
\small
\renewcommand{\arraystretch}{.8}
\begin{tabular}{p{.2\linewidth}|lp{.02\linewidth}}
\toprule
\textbf{Variables} & \textbf{Description} \\
\midrule
Population & $1{,}375$ seed population identifiers (Site 0), $1{,}194$ seed population identifiers (Site 1) \\
Original Planting Date & Actual planting date of the population \\
Early Planting Date & Earliest date the population can be planted \\
Late Planting Date & Latest date the population can be planted \\
Required GDUs & GDUs needed to harvest a seed population (range: 649 to $1{,}414$) \\
Harvest Quantity & Number of ears of corn for each seed population (Scenario 1 or Scenario 2) \\
Historical Daily GDUs & GDUs accumulated for each calendar day from 2009-2019 (Site 0 or Site 1) \\
\bottomrule
\end{tabular}
\caption{Description of the data provided in the 2021 Syngenta Crop Challenge.}
\label{tab:data_description}
\end{table}


As mentioned, the Challenge included two Scenarios, and for each Scenario participants were required to specify the planting date for each seed population. For Scenario 1, the weekly storage capacity was specified to be 7{,}000 ears at site 0 and 6{,}000 ears at site 1, with the sites to be optimized separately. The goal as specified in the Challenge was to minimize the first three quantities specified in Section \ref{sec:intro}, though as discussed we have added a fourth. For Scenario 2, the task was to specify the schedule but also provide a meaningful storage capacity.

\subsection{Mathematical Formulation} \label{sec:problem}

Suppose at site 0 there are $n_0$ seed populations and at site $1$ there are $n_{1}$ seed populations. Without loss of generality, we will describe the optimization setting for site $0$. We represent the $n_{0}$ seed populations as $s_{1},s_{2},\ldots, s_{n_{0}}$. Seed population $s_{i}$ is associated with a known harvest quantity $h(s_{i})$; a known amount of growing degree units (GDUs) necessary to trigger its harvest, $g(s_{i})$; and an interval composed of two dates within which the population must be planted, $[L(s_{i}),U(s_{i})]$. Also, let $g_{d}$ be the GDU for day $d$ for $d = 1,2,...,D$, where $D$ is the number of days in the optimization period. Furthermore, $\hat{g}_{d}$ is the number of predicted GDUs for day $d$, as estimated by our predictive model (see Section \ref{sec:predict}). Finally, we use $p_{i}$ to represent the planting day for seed population $s_{i}$. The $\mathbf{p}=p_{i}$ are the decision variables in our optimization and must be such that $L(s_{i}) \leq p_{i} \leq U(s_{i})$, for $i=1,2,\ldots,n_{0}$. We then have the full harvest quantity on week $j$ as $h_{j}(\mathbf{p}) = \sum_{i \in \mathcal{H}_{j}} h(s_{i}), \; j=1,2,\ldots W$, where $\mathcal{H}_{j}$ includes the indices for seed populations that are harvested in week $j$, $W$ is the number of weeks in the optimization period, which means $W$ should be counted from the first week with corn planted. Mathematically, and to emphasize $h_{j}$'s dependence upon $\mathbf{p}$, note that $i \in \mathcal{H}_{j}$ if $\sum_{d=p_{i}}^{d^{*}} \hat{g}_{d} > g(s_{i})$ where $d^* = \text{arg min}_{d'=p_{i},p_{i}+1,\ldots,D}\left\{ \sum_{d=p_{i}}^{d'} \hat{g}_{d} > g(s_{i})\right\}$. Finally, let $\mathcal{J}$ be the set of all week indices for which $h_{j}>0$ and $C_0=7{,}000$ be the weekly harvest capacity for site 0 in scenario 1 ($C_{1}=6{,}000$). 

\subsubsection{Model 1: Base} \label{sec:model1}

Given the problem specification above there are four criteria to minimize, and they can be represented as $\mathbf{f}_{1}(\mathbf{p})=(f_{11}(\mathbf{p}),f_{12}(\mathbf{p}),f_{13}(\mathbf{p}),f_{14}(\mathbf{p}))$, where
\begin{itemize}
    \item $f_{11}(\mathbf{p}) = \text{median}_{j \in \mathcal{J}}(|C_0-h_{j}(\mathbf{p})|)$, 
    \item $f_{12}(\mathbf{p})= \max_{j \in \mathcal{J}}(|C_0-h_{j}(\mathbf{p})|)$, 
    \item $f_{13}(\mathbf{p}) = \sum_{j=1}^{W} I(h_{j}(\mathbf{p})>0)$, 
    \item $f_{14}(\mathbf{p})= \sum_{j=1}^{W} (h_{j}(\mathbf{p})-C_{0})^{+}$,
\end{itemize}
$I$ is an indicator function, and $x^+=\max(x,0)$. 
The first and second criteria control the median and maximum absolute difference from capacity, respectively, while the third criterion minimizes the number of non-zero harvest weeks. The final criterion minimizes the total amount of wasted product. 

This basic criteria set has the possible drawback that the penalty for being overcapacity is not very severe, which may lead to difficulty in finding an efficient Pareto front based on this model. We address this issue in the soft constraint models in the next sections. 

\subsubsection{Model 2: First Penalty Model} \label{sec:model2}
 
For a given planting schedule, let $a_j=h_j-C_0, \; j \in \mathcal{J}_a$ where $\mathcal{J}_a$ is the set of weeks for which the harvest quantity is above capacity with cardinality $n_{a}$, and $b_j=C_0-h_j, \; j \in \mathcal{J}_b$ where $\mathcal{J}_b$ is the set of weeks for which the harvest quantity is below capacity with cardinality $n_{b}$. These quantities separately mark the positive part and negative part of the weekly harvest deviations from capacity. In order to impose a soft capacity constraint, we solve a new multiobjective optimization problem with the following criteria to be simultaneously minimized: $\mathbf{f}_{2}(\mathbf{p})=(f_{21}(\mathbf{p}),f_{22}(\mathbf{p}),f_{23}(\mathbf{p}),f_{24}(\mathbf{p}))$, where
\begin{itemize}
    \item $f_{21}(\mathbf{p}) = \text{median}_{j \in \mathcal{J}}|C_0-h_{j}(\mathbf{p})|$, 
    \item $f_{22}(\mathbf{p}) = \frac{\sum_{j\in \mathcal{J}_a}a^r_{j}(\mathbf{p})}{n_a}$,
    \item $f_{23}(\mathbf{p}) =  \frac{\sum_{j\in \mathcal{J}_b}b_{j}(\mathbf{p})}{n_b}$,
    \item $f_{24}(\mathbf{p}) = \sum_{j=1}^{W} I(h_{j}(\mathbf{p})>0)$
\end{itemize} 
As with the Base Model, this formulation controls both the median absolute deviation from capacity and the total number of nonzero harvest weeks. However, this model directly confronts the asymmetry between the costs of being overcapacity and undercapacity by minimizing the average of the $r^{\text{th}}$ power of all overcapacity amounts as well as the unexponentiated average of all undercapacity amounts.

\subsubsection{Model 3: Second Penalty Model} \label{sec:model3}

Finally, we consider an additional model that focuses on severely penalizing deviations from capacity, but also making the deviations as uniform as possible. In particular, Model 3 includes the following criteria, again to simultaneously minimize: $\mathbf{f}_{3}(\mathbf{p})=(f_{31}(\mathbf{p}),f_{32}(\mathbf{p}),f_{33}(\mathbf{p}),f_{34}(\mathbf{p}))$
\begin{itemize}
    \item $f_{31}(\mathbf{p}) = \text{median}_{j \in J}|C_0-h_{j}(\mathbf{p})|^{r}$, 
    \item $f_{32}(\mathbf{p}) = \text{max}_{j \in J}|C_0-h_{j}(\mathbf{p})|^{r}$,
    \item $f_{33}(\mathbf{p}) = \text{sd}_{j \in J}(C_0-h_{j}(\mathbf{p}))$,
    \item $f_{34}(\mathbf{p}) = \sum_{j=1}^{W} I(h_{j}(\mathbf{p})>0)$,
\end{itemize} 
where $\textbf{sd}$ represents the standard deviation of the harvest deviations from capacity. Note that we also make the median and max criteria more severe by adding a power to each.

\section{Solution Methodology} \label{sec:methods}

In Sections \ref{sec:predict} and \ref{sec:optimize} we provide a description of our prediction of GDUs and optimization approach, respectively. Overall, our solution strategy includes two basic building blocks. First, we construct a predictive model to predict the GDUs for each relevant day in 2020 and 2021, at each site. Then, we build a multiobjective optimization framework to construct Pareto fronts of solutions for the four objectives, using a genetic algorithm. In order to enrich our search of the solution space, we consider the three model formulations discussed in Section \ref{sec:problem}, and for each we perform an extensive round of tuning. This leads to an assessment of the resulting Pareto fronts, and we use the hypervolume indicator to choose between the three Pareto fronts (one for each of the three optimization models of Section \ref{sec:problem}). Finally, we select a particular solution from the chosen Pareto front by using the Technique for Order of Preference by Similarity to Ideal Solution \cite[TOPSIS,][]{topsiswang}. Note that when we compare the three Pareto fronts, we always do so by assessing their quality according to the Base Model (Section \ref{sec:model1}), since this is most reflective of the original Challenge objectives.

\FloatBarrier

\subsection{Predictive Models for Daily GDUs} \label{sec:predict}

A critical parameter within the optimization problem are the $g_{d}$'s, which are the Growing Degree Units (GDUs) for day $d$ within the planting window (January 1, 2020 to February 16, 2021) at site 0 and site 1. After seed population $s_i$ is planted, it will be harvested once the accumulated GDUs are greater than $g(s_i)$. Since we don't know the $g_d$'s in advance for both sites, we must estimate them from historical data, $\hat{g}_{d}$. As mentioned in Section \ref{sec:problem_data} and Table \ref{tab:data_description}, this data specifies the number of daily GDUs for each calendar day from 2009 to 2019 at site 0 and site 1. 
We compare two methods to predict the ${g}_{j}$'s. The first one is a simple averaging procedure, and the second is a long short-term memory (LSTM) neural network.

\subsubsection{Simple Averaging Model}

This method simply takes the average of the GDUs for each specific day of the year. For instance, in order to predict the daily accumulated GDUs on January 1, 2020, we take the average for all of January 1’s from 2009 to 2019. This average value is counted as the predicted value for January 1, 2020 and January 1, 2021. The same procedure is performed for the other 364 days of the year with all leap days being omitted.

\subsubsection{Long Short-Term Memory Model}

The second model we used for predicting daily GDUs is the Long Short-Term Memory (LSTM) neural network \citep{lstm}. As a type of deep learning method, LSTM can model both linear and non-linear relationships within time-dependent observations, which means it is more flexible than traditional time-series modeling methods such as ARIMA \citep{arimalstm,Box2013}. 
LSTM is an improved version of a recurrent neural network \citep[RNN,][]{medsker1999recurrent}, and the architecture of LSTM can address the limitations of traditional RNN on modelling long-term dependencies. In our work, we used the LSTM network to predict GDU $g_{ymd}$, where $y$ is the year, $m$ is the month, and $d$ is the day within month. (Note that our notation deviates slightly here, since previously as part of the problem specification we defined $g_{d}$ as the GDU for day $d$ of the optimization period; here, we use more flexible notation to capture days as far back as 2011.) We use this to define a sequence of GDU's comprised of the same day/month from the previous $l$ years. That is, $g_{ymd}$ is associated with the sequence $g_{y-l,md}$, \ldots, $g_{y-2,md}$, $g_{y-1,md}$. We denote this sequence $(g_{y-l,md}, \ldots, g_{y-2,md}, g_{y-1,md}; g_{ymd})$. We consider, as our training set, all such sequences for $y \in (2017, 2018, 2019)$ with $l=8$ (a total of $365 \times 3 = 1{,}095$ sequences). For instance, a particular sequence would be $(g_{2011,7,13}, g_{2012,7,13}, \ldots, g_{2018,7,13}; g_{2019,7,13})$ and this would use information from July 13 in years $2011-2018$ to predict the GDUs on July 13, 2019. We provide details of the procedure, along with the hyperparameters that we used, in the Supplementary material.

\subsubsection{Predictive Model Results}

\begin{figure}[hbt!]
    \centering
    \includegraphics[width=0.75\textwidth]{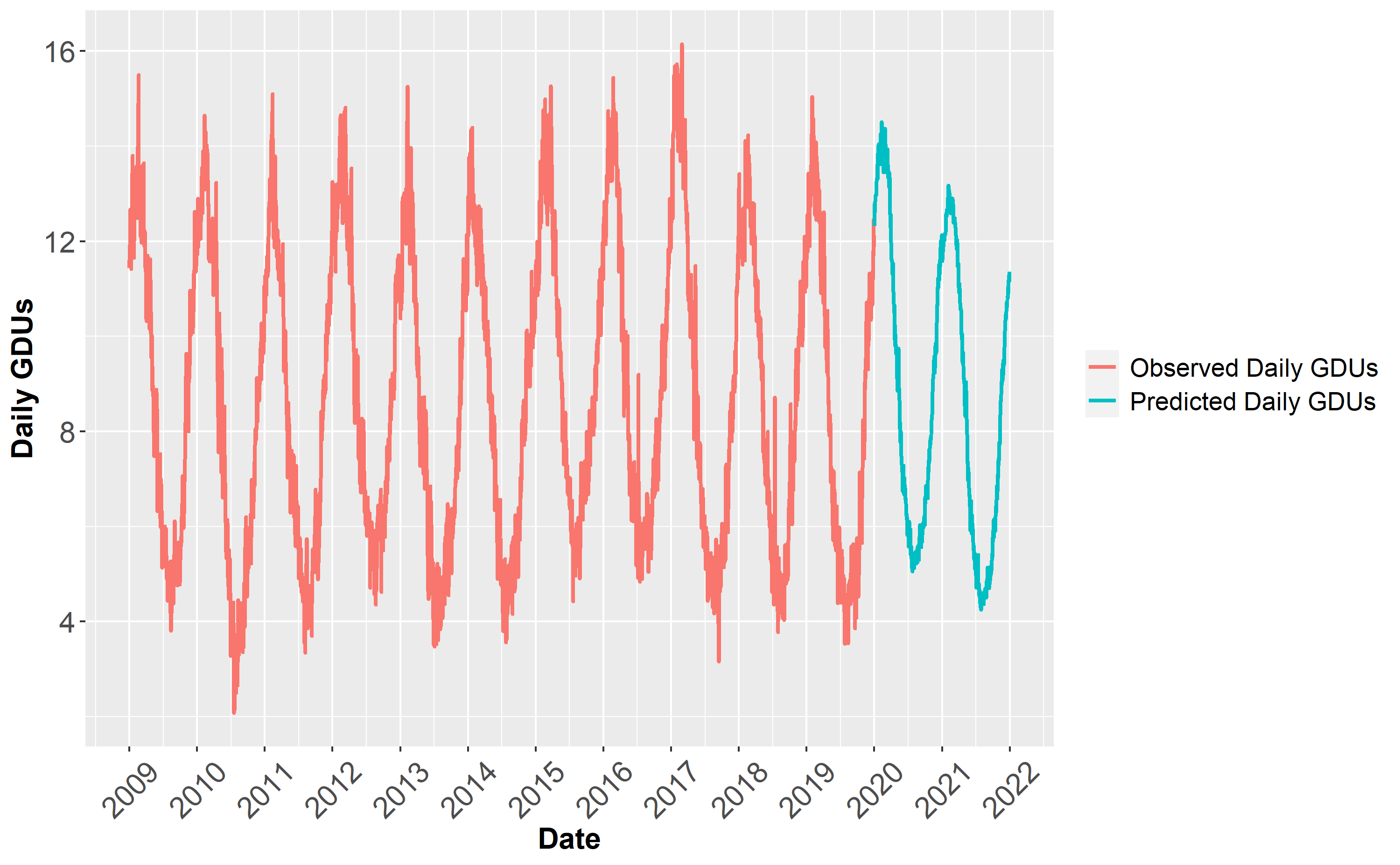}
    \caption{For Site 0, real GDU values between 2009 and 2019, and predicted GDU values for 2020 and 2021, using the LSTM(2,80) model.}
    \label{fig:0921PredGDU_site0}
\end{figure}

The models were trained and evaluated on a holdout set using the mean absolute error (MAE) and mean squared error (MSE).  Table \ref{table::lstmParam_site0} shows the performance of LSTM models for Site 0 and Site 1, respectively, computed by using the sequences implied by $y \in \{2017,2018\}$ as the training set, with the sequences in 2019 as the holdout set. Based on the Tables and the previous informal testing, we chose to use the LSTM model with 2 hidden layers, a batch size of 5, and 80 epochs. The final version of the model was trained using the sequences implied by $y \in \{2017,2018,2019\}$, and used to predict GDU's for 2020. When predicting time series, the further into the future predictions are made, the more variability is included in the predictions, and thus less seasonal structure is included in the predictions. In our case, instead of using the $y \in \{2017,2018,2019\}$-trained model to predict GDU's in 2021, we trained the model again with $y \in \{2018,2019,2020_{\hat{g}}\}$ and $\l=9$, where $2020_{\hat{g}}$ represents the predicted values for 2020. That is, we imputed GDU values for 2020 in order to make more variable predictions for 2021, so that the resulting predictions produced patterns more similar to past seasonal behavior. Figure \ref{fig:0921PredGDU_site0} visualizes the historical daily GDUs from 2009 to 2019 and the predicted daily GDUs from 2020 to 2021 for Site 0 (similar plot for Site 1 in Supplementary Material). 

\begin{table}[]
\begin{tabular}{l|ccc||ccc}
                     & LSTM(2, 80)    & LSTM(3, 80) & Simple Avg &  LSTM(2, 80)    & LSTM(3, 80) & Simple Avg \\ \hline
MAE & \textbf{0.548} & 0.759       & 0.624 &   \textbf{0.725} & 0.754       & 0.825         \\
MSE   & \textbf{0.682} & 0.911       & 0.757 &   \textbf{0.952} & 0.992       & 1.062 
\end{tabular}
\caption{Comparison of 2019 prediction results for Site 0 (left) and Site 1 (right), comparing LSTM Models and the simple averaging model. For LSTM(x, y) in the header, x is the number of layers of LSTM, and y is the number of epochs.}
\label{table::lstmParam_site0}
\end{table}


\FloatBarrier

\subsection{Multi-objective Optimization Methodology} \label{sec:optimize}

To solve each of the specified optimization problems in Section \ref{sec:problem}, we used an R implementation of the non-dominated sorting genetic algorithm II (NSGA-II) \citep{deb2002fast,nsga2R}. This procedure includes the elements of classical genetic algorithms \citep{whitley1994genetic,kumar2010genetic,mirjalili2019genetic,katoch2021review}, including the characterization of solutions as genes, as well as crossover and mutation. 
NSGA-II naturally handles multiple objectives in settings where objectives are treated as a black box, and attempts to provide a set of Pareto optimal solutions. 
In this section we provide a description of our solution strategy, including an extensive search of the tuning parameter space, details regarding Pareto optimality, the hypervolume measure of Pareto fronts, and the TOPSIS approach to choosing a final solution.

\subsubsection{Solution Strategy} \label{sec:strategy}

We characterize the solution of our multiobjective optimization problem in terms of Pareto optimality. Following \citet{cao2015using} and using notation from Sections \ref{sec:model1}-\ref{sec:model3}, for optimization model $\mathbf{f}_{k}$, a solution $\mathbf{p}_{1}$ is said to dominate solution $\mathbf{p}_{2}$ if $f_{ko}(\mathbf{p}_{1}) \leq f_{ko}(\mathbf{p}_{2}) \; \forall \; o\in \{1,2,3,4\}$ with $f_{ko}(\mathbf{p}_{1}) < f_{ko}(\mathbf{p}_{2})$ for at least one $o \in \{1,2,3,4\}$; whereas $\mathbf{p}_{1}$ weakly dominates $\mathbf{p}_{2}$ if $f_{ko}(\mathbf{p}_{1}) \leq f_{ko}(\mathbf{p}_{2}) \; \forall \; o \in \{1,2,3,4\}$. A set of solutions is Pareto optimal if the set consists of nondominated solutions, and this set of solutions is said to occupy the Pareto front. The hypervolume measure \citep{zitzler1998multiobjective} measures the volume of the criteria-space that is weakly dominated by the set of points composing a Pareto front, and we use this measure to compare the quality of different Pareto fronts. In order to bound the hypervolume measure, a reference point must be defined and we discuss that below.


One difficulty with the NSGA-II approach is that a number of hyperparameters must be specified in order to ensure an effective solution: the number of generations, the population size, and both crossover and mutation probabilities. In addition, for the two penalty models (Sections \ref{sec:model2} and \ref{sec:model3}), we denote the exponent $r$ as an additional parameter to be chosen. In the following subsection we provide more details regarding our study of these parameters, which explores a variety of parameter combinations for each of the three multi-objective problems. For now we simply denote the set of parameter combinations---for each of $\mathbf{f}_{1}$, $\mathbf{f}_{2}$, and $\mathbf{f}_{3}$---as $\mathcal{P}_{1}$, $\mathcal{P}_{2}$, and $\mathcal{P}_{3}$, respectively. Each implementation of NSGA-II for $\mathbf{f}_{k}$, for each parameter combination in $\mathcal{P}_{k}$, will include a final population of solutions, denoted $\text{POP}_{k\ell}, \; k=1,2,3; \; \ell=1,2,...,c_{k}$, where $c_{k}$ is the number of parameter combinations for $\mathcal{P}_{k}$. Each population $\text{POP}_{k\ell}$ includes $n_{k\ell}$ solutions, depending on the population size of the NSGA-II implementation for model $k$ and parameter combination $\ell$. 
Our basic solution strategy, then, is:
\begin{enumerate}
    \item Solve $\mathbf{f}_{k}$ for each combination of tuning parameters in $\mathcal{P}_{k}$, $k=1,2,3$, respectively, yielding $\text{POP}_{k\ell}, \; k=1,2,3; \; \ell=1,2,...,c_{k}$.
    \item Construct three populations $P_k = \bigcup\{ \text{POP}_{k\ell} |  \ell \in \{1,2,\ldots,c_{k}\}\}$ with $k = 1, 2, 3$. Map each $\mathbf{p} \in P_k$ to the criteria in $\mathbf{f}_{1}$ as $\mathbf{f}_{1}(\mathbf{p})$, because we are evaluating all solutions, regardless of whether they solved $\mathbf{f}_{1}$, $\mathbf{f}_{2}$, or $\mathbf{f}_{3}$, in terms of the objective functions associated with $\mathbf{f}_{1}$. 
    \item From each $P_k$ after mapping to the criteria $f_1$, construct Pareto fronts $PF_k, k = 1, 2, 3$. Then construct a meta-population $P_{PF} = \bigcup\{ PF_k | k \in \{1,2,3\}\}$. 
    \item Let $F^{\text{min}}_{1o}=\min_{\mathbf{p} \in P_{PF}}f_{1o}(\mathbf{p})$ and $F^{\text{max}}_{1o}=\max_{\mathbf{p} \in P_{PF}}f_{1o}(\mathbf{p})$ be the minimum and maximum value of criterion $f_{1o}$, respectively, for model $1$ and objective $o$, across all solutions in $P_{PF}$. For each $\mathbf{p} \in P_{PF}$, scale $\mathbf{f}_{1}(\mathbf{p})$ such that $\mathbf{f}_{1}^{\text{sca}}(\mathbf{p}) \in [0,1]^{4}$, where $f^{\text{sca}}_{1o} = \frac{f_{1o}-F^{\text{min}}_{1o}}{F^{\text{max}}_{1o}-F^{\text{min}}_{1o}}$ for $o=1,2,3,4$.
    \item For each $\text{PF}_{k}$ after scaling in Step 4, compute hypervolume $v_{k}$ based on $[\mathbf{f}_{11}^{sca}, \mathbf{f}_{12}^{sca}, \mathbf{f}_{13}^{sca}, \mathbf{f}_{14}^{sca}]$, and choose $k^* = \text{argmax}_{k} v_{k}$. This defines the model which yields the best scaled Pareto front.
    \item From $\text{PF}_{k^*}$, the best single solution $\mathbf{p}^{*}$ is chosen via TOPSIS \cite[][see Section \ref{sec:TOPSIS}]{topsiswang}. The chosen solution $\mathbf{p}^*$ produces $\mathbf{f}_{1}(\mathbf{p}^*)$ as the final criteria values from the originally specified set of objectives.
\end{enumerate}

\textbf{Remark 1}: We emphasize that Models 2 and 3 ($\mathbf{f}_{2}(\mathbf{p})$ and $\mathbf{f}_{3}(\mathbf{p})$, respectively) are simply vehicles to improve the solutions obtained with respect to Model 1, because Model 1 most closely aligns with the original requirements of the Competition. This is why, in Step 2 and Step 6 above, we evaluate the solutions using the objectives defined in Model 1.

\textbf{Remark 2}: In order to compute the hypervolumes $v_{k}$, we not only need to scale the criteria, but we also must define a reference point. After some experimentation, we chose $[2,2,2,2]$, recognizing that for all of the reference points we explored, the hypervolume ordering of the Pareto fronts stayed the same.

\textbf{Remark 3}: We have chosen to take a middle road regarding the pooling of solutions. On the one hand, one could choose to pool all solutions from all models, and construct a single Pareto front from which a final solution is chosen, so that the hypervolume technique can be ignored. On the other hand, one could form separate Pareto fronts for each parameter combination for each model, so that the best Pareto front is chosen based on both the best model and parameter combination. Instead, we have pooled across parameter combinations but preserving separate populations of solutions for each model. This allows some insight to be obtained regarding whether the alternative models result in more attractive solutions. 

\FloatBarrier

\subsubsection{Tuning Parameter Selection}

The NSGA-II heuristic requires the input of several parameters: the number of generations, the population size, and both crossover and mutation probabilities. In addition, for the two penalty models, we denote the exponent $r$ as an additional parameter. In order to study the quality of the Pareto fronts across values of these parameters, we used NSGA-II for a number of different combinations of these parameter values, denoted $\mathcal{P}_{k}$ for model $k$ above. \autoref{tab:tuning_pars} shows these sets, chosen based upon several works in the literature \citep{alender, boyabtli, gaparamtuning}. \cite{alender} and \cite{gaparamtuning} suggested exploring from $N$ to $3N$, where $N$ is the problem size (in our case, the total number of seed populations). So, possible  population sizes for site 0 are $N=1{,}376$, $2N=2{,}752$ and $3N=4{,}128$, and similarly for site 1. For the other settings of tuning parameters, we are guided in part by the design of parameters in \cite{gaparamtuning}. In all, $\mathcal{P}_{1}$ consisted of $3^4=81$ parameter combinations, while $\mathcal{P}_{2}$ and $\mathcal{P}_{3}$ included $3^5=243$. This was a computationally intensive process, so we used Miami University's Redhawk Cluster, which has a CentOS 7.9.2009 operating system and uses an Intel Xeon Gold 6126 processor, with 2.6 GHz clock speed. Each node of the cluster has 96 GB available memory and the optimization problem for each parameter combination was solved using a single CPU core.

\begin{table}[]
\centering
\caption{For site 0, the three levels studied for each tuning parameter.}
\begin{threeparttable}
\begin{tabular}{llll}
  & \textbf{Low} & \textbf{Middle} & \textbf{High} \\ \hline
\textbf{Crossover Rate} & 0.5 & 0.75 & 1.0 \\
\textbf{Mutation Rate} & 0.001 & 0.01 & 0.1 \\
\textbf{Population Size}\tnote{1} & 1376 & 2752 & 4128 \\
\textbf{Generation Size} & 8000 & 10000 & 12000 \\
\textbf{Penalty Power}\tnote{2} & 1 & 2 & 3 \\ \hline    
\end{tabular}
\begin{tablenotes}
    \item [1] Only for Models 2 and 3.
    \item [2] For site 1, population sizes were $1{,}196$, $2{,}392$, and $3{,}588$.
\end{tablenotes}
\end{threeparttable}
\label{tab:tuning_pars}
\end{table}

Thus, we use NSGA-II to optimize each model for each of the tuning parameter combinations in $\mathcal{P}_{1}$, $\mathcal{P}_{2}$ and $\mathcal{P}_{3}$, respectively. To determine which model, and which set of tuning parameter values, produced the best Pareto front, we implemented the strategy outlined in Section \ref{sec:strategy}. Note that we tried several approaches to determine which parameter combination would optimize the hypervolume, include response surface methodology \citep{myers2016response};
however, these methods did not clearly improve upon the results observed directly by choosing the parameter combination, for each model, that produced the best hypervolume.

\FloatBarrier

\subsubsection{Choosing the Final Solution from Pareto Front} \label{sec:TOPSIS}

The Technique for Order of Preference by Similarity to Ideal Solution \citep[TOPSIS,][]{topsiswang} is a straightforward method to choose a particular Pareto solution. In Scenario 1, TOPSIS chooses the final solution $\mathbf{p}^*$ from the Pareto front $PF_{k^{*}}$. For simplicity, in this section we will explain the procedure using a generic Pareto front $PF$, composed of $i=1,2,\ldots,n$ solutions each evaluated on $o=1,2,3,4$ objectives. The basic idea of the method: TOPSIS first normalizes each solution in the Pareto front (Step 1 below), and then uses the Euclidean distance of a solution from both the positive and negative ideal to choose a particular solution (Steps 3-5 below).

\begin{itemize}
  \item Step 1. For solution $i$ and objective $o$, normalize the Pareto front to create the normalized Pareto front $F_{io}$ as follows:
  
  $$F_{io} = \frac{PF_{io}}{\sqrt{\sum_{i=1}^{n}PF_{io}^2}}$$
  
  \item Step 2. Adjust each solution according to user-specified weights $w_{o}$: 
  $$V_{io} = F_{io} * w_{o}$$
  For Scenario 1, we assume each objective function is equally important, thus the weights are constant.
  
  \item Step 3. Find the positive ideal solution $A^+$ and the negative ideal solution $A^-$. Here, our optimization goal is to minimize each objective; i.e. $min(\mathbf{f}_{1o}(p)) \  \forall o \in \{1,2,3,4\}$, so that $A^+ = \{ min(V_{io}) | o \in \{1,2,3,4\} \}$ and $A^- = \{ max(v_{io}) | o \in \{1,2,3,4\} \} $.
  
  \item Step 4. Compute the Euclidean distance $S_{i+}$ from each point $V_{i}$ to the positive ideal solution $A^+$, and the Euclidean distance $S_{i-}$ from each point $V_{i}$ to the negative ideal solution $A^-$.

  \item Step 5. Calculate the score $C_i = \frac{S_{i-}}{S_{i-} + S_{i+}}$ and choose the solution with the largest $C_i$ as the final recommended solution.
  
\end{itemize}
As can be seen, solutions close to the positive ideal and far from the negative ideal will be chosen.

There are several other approaches to choosing a final solution that might be used here, including approaches suggested by \citet{LINMAP}, \citet{vikor}, \citet{saw}, and \citet{topsiswang}. \cite{topsiswang} suggest that TOPSIS is the most commonly used method for choosing the final single optimal solution of a multi-objective optimization problem. We chose it because it is relatively simple and doesn't require many user inputs.



\FloatBarrier
\subsection{Scenario 2 Description and Methodology} \label{sec:scenario2_description}

In contrast to Scenario 1, Scenario 2 requires a generalization of our methods in order to estimate and optimize the location capacity---a quantity we will denote $\hat{C}$---instead of assuming a fixed capacity, $C_{0}$, all while optimizing the planting schedule as well. While the predictive modeling approach for daily cumulative GDUs will stay the same (see Section \ref{sec:predict}), we generalize the optimization by treating both the location capacity and planting schedule as decision variables, while adding the location capacity  as an objective function as well. Due to space constraints, we provide the details of our formulation in the Supplementary Material.

\FloatBarrier
\section{Results} \label{sec:results}

In this section, we provide solutions for Scenario 1, Site 0 using the methodology provided in Section \ref{sec:optimize}. We provide results for Scenario 1, Site 1 and Scenario 2, Site 0 in the Supplementary material. We could use the same methodology to provide results for Scenario 2, Site 1, but we omit them due to the large Pareto fronts and associated computational cost.

\begin{figure}[hbt!]
    \centering
         \centering
         \includegraphics[width=\textwidth]{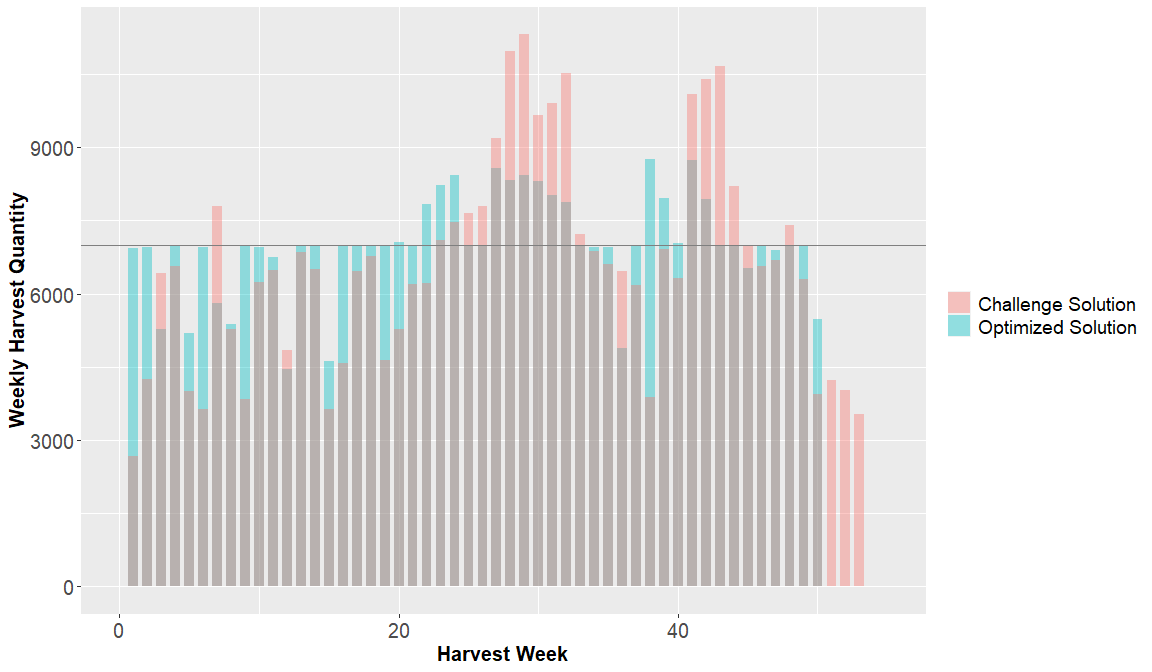}
         \caption{Comparison between the solution submitted to the Challenge and the solution optimized using the methods described in this paper.}
     \label{fig:Scenario1_Site0}
\end{figure}

In Scenario 1, a storage capacity is provided, and for Site 0 it is 7{,}000 ears per week. For this initial set of results, we provide detailed results in the Supplementary material. Table \ref{tab:result_s1_site0} and Figure \ref{fig:Scenario1_Site0} show our final solution, which Pareto-dominates the solution we originally submitted to the Challenge. Notice that our preferred solution includes three fewer harvest weeks, is clearly less variable around the capacity, and results in less waste. We have provided results for Scenario 1 (Site 1), along with Scenario 2 (Site 0) in the Supplementary Materials.

\begin{table}[h]
\centering
\caption{Criterion values for solutions $\mathbf{p}^*$ (final single optimal solution), $\mathbf{p}_{challenge}$ (solution submitted to the Challenge), $\mathbf{p}_{initial}$ (initial solution given by the Challenge), under Scenario 1, Site 0}
\label{tab:result_s1_site0}
\begin{tabular}{llll}
  & \textbf{$\mathbf{p}^*$} & \textbf{$\mathbf{p}_{challenge}$} & \textbf{$\mathbf{p}_{initial}$} \\ \hline
\textbf{Median Absolute Difference ($\mathbf{f}_{11}$)} & 53 & 823 & 4{,}012 \\
\textbf{Max Absolute Difference ($\mathbf{f}_{12}$)} & 2{,}537 & 4{,}320 & 19{,}222 \\
\textbf{\# of Non-zero Harvest Week ($\mathbf{f}_{13}$)} & 50 & 53 & 54 \\
\textbf{Total Amount of Wasted Product ($\mathbf{f}_{14}$)} & 16{,}514 & 34{,}320 & 102{,}204 \\ 
\end{tabular} 
\end{table}

\FloatBarrier

\section{Discussion and Conclusion} \label{sec:conclusion}

In this paper, we have provided a solution to a challenging corn-planting problem presented in the 2021 Syngenta Crop Challenge. In particular, we have improved upon the solutions submitted to the Challenge, based upon the three criteria provided along with an additional ``waste'' criterion. This problem requires multiobjective optimization and predictive modeling, and we improved the results by considering two related optimization problems that enhance solutions with respect to the original four criteria. We compare results across these three models using the hypervolume indicator, and choose a final solution using the TOPSIS method from the literature. In the end, we obtain solutions that clearly improve upon those our team submitted to the Challenge. (Note that the Challenge entry was submitted by the first three authors.) Compared to the Challenge entry, we extensively explored the tuning parameter space for the Genetic Algorithm, added an additional objective function to the base model ($\mathbf{f}_1$), modified the objective functions in Model 2 ($\mathbf{f}_2$), and added a new model ($\mathbf{f}_3$). In our original entry, we also chose a final solution using a weighted sum of normalized objectives, rather than the hypervolume measure and TOPSIS that we use here. 

As far as we know, this problem  has only been considered by two other Challenge teams \citep{sajid2022optimizing, lizhiwang}. \cite{sajid2022optimizing} used a Convolutional Neural Network to predict GDUs, and developed a Mixed Integer Linear Programming-based algorithm to determine a final solution. Like us, \cite{lizhiwang} uses LSTM for predictive modeling, whereas instead of optimizing the objectives individually, they focus on a single objective which minimizes the sum of the absolute difference between weekly harvest quantity and location capacity, reducing the problem to a single objective optimization problem with other objectives ignored. Table \ref{tab:res_comparison} shows the comparison of criteria values between our final optimized solution $\mathbf{f}_{1}(\mathbf{p}^*)$ and the optimized solutions made by \cite{sajid2022optimizing, lizhiwang} under Scenario 1 as well as Scenario 2, Site 0. No group's solutions uniformly dominate any other, but for Scenario 1 especially our solutions are strong. We do uniformly outperform the other teams in the Number of Harvest Weeks criterion ($\mathbf{f}_{13}$).

There are several final points to make regarding this work. First, it represents a full solution to a challenging, real problem, which includes predictive modeling and multiobjective optimization. Note that we handled an inherently multiobjective optimization problem with multiobjective optimization tools, instead of combining into one objective or solving the various objectives separately. Instead, we used a genetic algorithm that explicitly handles the trade-off between objectives while considering them simultaneously. Secondly, we demonstrate the complexity inherent in employing multiobjective optimization in such a problem. Not only did we undertake substantial exploration of the tuning parameter space, we also proposed a procedure by which Pareto fronts were compared, and a final solution chosen. We also demonstrate that in this case, better solutions to the optimization problem of interest (Model 1) are obtained by estimating the solution to a related but distinct problem (Model 2). 

\begin{table}[] 
\caption{Comparison of final solutions between our proposed solution ($\mathbf{f}_{1}(\mathbf{p}^{*})$) and the published results of two other teams.}
\label{tab:res_comparison}
\resizebox{\textwidth}{!}{%
\begin{tabular}{|l|l|l|l|l|l|}
\hline
Scenario and Site                   & Solution Source                 & $\mathbf{f}_{11}(\mathbf{p^*})$ & $\mathbf{f}_{12}(\mathbf{p^*})$ & $\mathbf{f}_{13}(\mathbf{p^*})$ & Estimated Capacity \\ \hline
\multirow{3}{*}{Scenario 1, Site 0} & $\mathbf{f}_{1}(\mathbf{p^*})$                       & 53       & 2{,}537     & 50       & -                  \\ \cline{2-6} 
                                    & \cite{lizhiwang} & 57       & 2{,}471     & 51       & -                  \\ \cline{2-6} 
                                    & \cite{sajid2022optimizing} & 16.5     & 2{,}883     & 51       & -                  \\ \hline
\multirow{3}{*}{Scenario 1, Site 1} & $\mathbf{f}_{1}(\mathbf{p^*})$                         & 8        & 827      & 51       & -                  \\ \cline{2-6} 
                                    & \cite{lizhiwang}  & 25       & 2{,}633     & 52       & -                  \\ \cline{2-6} 
                                    & \cite{sajid2022optimizing} & 29.5     & 394      & 52       & -                  \\ \hline
\multirow{3}{*}{Scenario 2, Site 0} & $\mathbf{f}_{1}(\mathbf{p^*})$                         & 409.6  & 4{,}286.6 & 51       & 10{,}312.6            \\ \cline{2-6} 
                                    & \cite{lizhiwang}  & NA       & NA       & NA       & 10{,}795              \\ \cline{2-6} 
                                    & \cite{sajid2022optimizing} & 49.5     & 2{,}057     & 52       & 9{,}800               \\ \hline
\end{tabular}%
} 
\end{table}

\FloatBarrier




\section*{Acknowledgements}

The authors are grateful to Syngenta for providing the opportunity to work on this problem, as well as to the committee who judged the initial submissions to the Challenge. 

\section*{Conflicts of Interest}

The authors have no funding support to acknowledge nor competing interests to declare.

\bibliography{mybibfile}

\end{document}